\documentclass [11pt, a4paper, twoside]{article}

\usepackage {diagrams, a4, indentfirst, latexsym, theorem, QED, amssymb, epsfig}

\renewcommand{\qedsymbol}{\quad\Box}
\renewcommand{\TheWordProof}{\bf Proof.}
\newcommand\nothing[1]{\relax}

\theoremstyle{change}
\theorembodyfont{\it}
\newtheorem{theorem}{Theorem.}[section]
\newtheorem{proposition}[theorem]{Proposition.}
\newtheorem{lemma}[theorem]{Lemma.}
\newtheorem{corollary}[theorem]{Corollary.}
{\theorembodyfont{\upshape}
\newtheorem{definition}[theorem]{Definition.}

\newtheorem{example}[theorem]{Example.}

\newtheorem{numpar}[theorem]{}

}

\newcommand{\bR}{{\mathbb R}}
\newcommand{\bZ}{{\mathbb Z}}
\newcommand{\bC}{{\mathbb C}}
\newcommand{\iso}{\cong}
\newcommand{\tensor}{\otimes}
\newcommand{\id}{\mathrm{id}}
\newcommand{\ie}{{\it i.e.}}

\newcommand{\fan}{\Sigma}

\date{}

\begin{document}

\title{A cohomological interpretation of Brion's formula}
\author{Thomas H\"uttemann}
\maketitle

\centerline {\it University of Leicester, Department of Mathematics}
\centerline {\it University Road, Leicester LE1~7RH, England (UK)}
\centerline {e-mail: \texttt{th68@mcs.le.ac.uk}}

\vglue 2\bigskipamount \hrule \medskip

{\footnotesize \smallskip

  A subset~$K$ of $\bR^n$ gives rise to a formal \textsc{Laurent\/}
  series with monomials corresponding to lattice points in~$K$. Under
  suitable hypotheses, this series represents a rational
  function~$R(K)$. \textsc{Michel Brion} has discovered a surprising
  formula \cite{Brion-latticepoints} relating the rational function
  $R(P)$ of a lattice polytope~$P$ to the sum of rational functions
  corresponding to the supporting cones subtended at the vertices
  of~$P$. The result is re-phrased and generalised in the language of
  cohomology of line bundles on complete toric varieties.
  \textsc{Brion}'s formula is the special case of an ample line bundle
  on a projective toric variety.---The paper also contains some
  general remarks on the cohomology of torus-equivariant line bundles
  on complete toric varieties, valid over arbitrary commutative ground rings.


  \noindent {\it Additional keywords:\/} Polytope, cone, lattice point
  enumerator, toric variety, line bundle, \v Cech
  cohomology \hfill (\today)}

\medskip
\hrule
\vglue 2\bigskipamount

\section{Introduction}

The main result of this paper is a generalisation of a formula discovered by
\textsc{Brion}, relating the lattice point enumerator of a rational polytope
to the lattice points enumerators of supporting cones subtended at its
vertices \cite[\S2.2]{Brion-latticepoints} (see \cite{Beck:Brion_formulae} for
an introduction to the theory, and \cite{H-Brion} for an elementary geometric
proof). In spirit the proof of the generalisation is similar to
\textsc{Brion}'s original exposition, but avoids the use of equivariant
$K$-theory in favour of a more elementary treatment of cohomology of line
bundles on complete toric varieties.

Since line bundles are encoded by support functions defined on a fan, the
result can be re-formulated in combinatorial terms. This has been done for
upper convex support functions (corresponding to line bundles which are
generated by global sections) by \textsc{Ishida}
\cite[Theorem~2.3]{Ishida-polyhedral}, generalising the original result of
\textsc{Brion}. The present paper goes one step further and includes the case
of arbitrary, non-convex support functions.

We will give a precise formulation of the result below. Roughly
speaking, we prove that a sum of certain rational functions, all given
by infinite \textsc{Laurent} series, degenerates to a \textsc{Laurent}
polynomial, and interpret the coefficients of the occurring monomials
as homogeneous \textsc{Euler} characteristics of the sheaf cohomology
of an algebraic line bundle.

The proof relies on a non-standard computation of the cohomology of line
bundles on complete toric varieties (Theorem~\ref{thm:line_bundle_cohomology})
which is similar to, but slightly easier than, the standard result as given by
\textsc{Oda} \cite[Theorem~2.6]{Oda-Toric}. This computation in turn depends
on a variant of \textsc{\v Cech} cohomology
(Proposition~\ref{prop:CechCohomology}) which should be well-known; since it
seems not to be well-documented in available publications, we include a proof
at the end of the paper.

\subsubsection*{Notational conventions and the main result}

We have to introduce some notation first. Let $M \iso \bZ^n$ be a
lattice of rank~$n$. We call the set of maps $S = \mathrm{map} (M, \bC)$ the
set of {\it formal \textsc{Laurent\/} series}. Given an element $\mathbf{b}
\in M$ we let $x^\mathbf{b} \in S$ denote the map which is zero on~$M
\setminus \{\mathbf{b}\}$, and takes the value~$1$ on~$\mathbf{b}$. We
call~$x^\mathbf{b}$ the {\it \textsc{Laurent\/} monomial} with
exponent~$\mathbf{b}$.

The terminology can be justified.  Given a choice of basis $e_1, e_2, \ldots,
e_n$ of~$M$ we can write every element $\mathbf{b} \in M$ uniquely as
$\mathbf{b} = \sum_j b_j e_j$ with $b_j \in \bZ$. Then for $f \in S$
the formal sum
\[\sum_{\mathbf{b} \in M} f(\mathbf{b}) \cdot x_1^{b_1} x_2^{b_2}
\ldots x_n^{b_n}\]
is a \textsc{Laurent} series in the indeterminates $x_1, x_2, \ldots, x_n$.
The map $x^\mathbf{b}$ corresponds to the product $x_1^{b_1} x_2^{b_2}
\ldots x_n^{b_n}$, \ie, a series with a single non-trivial summand.

Let $P \subset S$ denote the subset of maps with finite support; in
particular, it contains the maps $x^\mathbf{b}$ defined above. After choosing
a basis of~$M$ we can identify $P$ with the ring of \textsc{Laurent\/}
polynomials in $n$ indeterminates; on the level of maps, the product is given
by a convolution formula. The same formula equips $S$ with the usual structure
of a $P$-module.

Set $M_\bR = M \tensor \bR \iso \bR^n$. We consider $M$ as a subset
of~$M_\bR$ using the natural identification $M = M \tensor 1$.
Given a subset $K \subseteq M_\bR$
and an element $b \in M_\bR$ we define
\[b+K = \{ b + x \,|\, x \in K \} \qquad \hbox{and} \qquad -K = \{ -x \,|\, x
\in K \} \ .\]

\begin{definition}
  For a subset $K \subseteq M_\bR$ we define the formal \textsc{Laurent}
  series
  \[R[K] = \sum_{\mathbf{a} \in M \cap K} x^{\mathbf{a}}
  \ \in S \ .\]
\end{definition}

\noindent A straightforward calculation shows $R[\mathbf{b}+K] =
x^\mathbf{b} R[K]$ for any $\mathbf{b} \in M$.

\medbreak

In favourable cases, for example when $K$ is a pointed rational
polyhedral cone in~$M_\bR$, the series $R[K]$ represents a rational
function (an element in the quotient field~$Q(P)$ of~$P$) which we
will denote $R(K) \in Q(P)$.

As an explicit example, for $K = \bR_{\leq 2} = 2 + \bR_{\leq 0} \subset \bR$,
we have $R[K] = x^2 R[\bR_{\leq 0}] = x^2 \sum_{a \leq 0} x^a$, so $R(K) =
x^2/(1-x^{-1})$. See \cite{Beck:Brion_formulae} for more examples.

\medbreak

Let $N = \hom_\bZ (M, \bZ) \iso \bZ^n$ be the dual lattice of~$M$.
Then $N_\bR = N \tensor \bR \iso \bR^n$ is naturally the dual of the
$\bR$-vector space $M_\bR$. The duals of $N$ and $N_\bR$ are
canonically isomorphic to $M$ and $M_\bR$, respectively.

Let $\fan$ be a finite complete fan in~$N_\bR$, consisting of strongly convex
rational polyhedral cones, and denote by $X_\fan$ the associated toric variety
defined over~$\bC$. (See \cite{Oda-Toric} for details on cones, fans, and the
relation to varieties.) Let $h \colon N_\bR \rTo \bR$ be a support function
on~$\fan$; on each cone $\sigma \in \fan$ it coincides with a linear function
$h_\sigma \tensor \mathrm{id}_\bR$ for some $h_\sigma \in \hom_\bZ (N,
\bZ) = M$. Define the rational function
\[R(\fan, h) = \sum_{{\sigma \in \fan} \atop {\dim \sigma = n}}
R(-h_\sigma + \sigma^\vee)\]
where $\sigma^\vee = \{ x \in M_\bR \,|\, \forall y \in \sigma:
\langle x,y \rangle \geq 0\}$ is the dual cone, defined using the
standard evaluation pairing $\langle x , y \rangle = y(x)$.
Denote the (algebraic) line bundle on~$X_\fan$ associated
to~$h$ by~$L_h$; see \ref{par:L_h} below for an explicit description.
Given a vector $\mathbf{a} \in M$ we write $H^k (X_\fan; L_h)_\mathbf{a}$ for the
homogeneous part of degree~$\mathbf{a}$ of the $k$th sheaf cohomology
of~$L_h$; see~\ref{par:hom_piece_cohomology}--\ref{par:desc_hom_chain_cx} below for an elementary
description.

\begin{theorem}
  \label{thm:generalised_Brion}
  The rational function $R(\fan,h)$ is a \textsc{Laurent} polynomial.  The
  coefficient of the monomial $x^\mathbf{a}$ in the polynomial $R(\fan,h)$ is
  the \textsc{Euler} characteristic of $H^* (X_\fan; L_h)_\mathbf{a}$, so it is
  given by the alternating sum
  \[\chi (L_h)_\mathbf{a} = \sum_{k=0}^n (-1)^k \dim_\bC\, H^k
  (X_\fan; L_h)_\mathbf{a} \ .\]
  In short, we have the equality
  \begin{equation}
    \label{eq:main_res}
    R(\fan, h) = \sum_{\mathbf{a} \in M} \chi (L_h)_\mathbf{a} \cdot
    x^\mathbf{a} \ .
  \end{equation}
\end{theorem}

\subsubsection*{Brion's formula for lattice polytopes}

For an $n$-dimensional polytope $K \subset M_\bR$ with vertices
in~$M$, let $\fan_K$ denote the inner normal 
fan of~$K$. The support function of~$K$, given by
\[h_K \colon N_\bR \rTo \bR, \quad x \mapsto - \inf \{ \langle p,x\rangle \,|\, p \in K\}\ ,\]
defines a support function on~$\fan_K$. Explicit calculation of sheaf
cohomology shows
\[H^k (X_{\fan_K}; L_{h_K}) = \cases{0 & if \( k \not= 0\) \cr
\bigoplus_{M \cap K} \bC & if \(k=0\)}\]
(see \cite[Corollary~2.9]{Oda-Toric}, or
\cite[Theorem~2.5.3]{H-nonlin_toric} for an elementary proof).
If $\sigma$ is an $n$-dimensional cone in~$\fan_K$, and $h_\sigma$ is the linear
function associated to~$h_K$ and~$\sigma$, then $-h_\sigma + \sigma^\vee$ is
the support cone of~$K$ subtended at the vertex corresponding to~$\sigma$.
Thus Theorem~\ref{thm:generalised_Brion} reduces to the original theorem
of \textsc{Brion} \cite[\S2.2]{Brion-latticepoints}: The
rational function $R(\fan_K, h_K)$ is a \textsc{Laurent} polynomial with
terms corresponding to the integral points of~$K$.

Similarly, by considering the support function $-h$, and using the
calculation
\[H^k (X_{\fan_K}; L_{-h_K}) = \cases{0 & if \( k \not= n\) \cr
\bigoplus_{M \cap \mathrm{int}\, (-K)} \bC & if \(k=n\)}\]
(see \cite[Corollary~2.8]{Oda-Toric}, or
\cite[Theorem~2.5.3]{H-nonlin_toric} for an elementary proof),
we see that the rational function $R(\fan_K, -h_K)$ is a \textsc{Laurent}
polynomial with summands corresponding to the integral points in the
interior of~$-K$, up to a factor of $(-1)^n$
\cite[\S2.5]{Brion-latticepoints}.

Finally, by considering a globally linear support function $h = \mathbf{a} \in
M$, so that
$L_\mathbf{a} \iso \mathcal{O}_{X_{\fan_K}}$, Theorem~\ref{thm:generalised_Brion} together
with the calculation
\[H^k (X_{\fan_K}; L_\mathbf{a}) = \cases{0 & if \( k
  \not= 0\) \cr \bC & if \(k=0\)} \]
with $H^0$ concentrated in homogeneous degree~$\mathbf{a}$ (see
\cite[Theorem~2.5.3]{H-nonlin_toric} for the case $\mathbf{a}=0$, the general
case follows easily) yield the identity $R(\fan_K, \mathbf{a}) = x^\mathbf{a}$
(cf.~\cite[Corollary~2.4]{Ishida-polyhedral}; see
\cite[Proposition~3.1]{Brion-Vergne} for a \textsc{Laurent} series version).

\medbreak

The cohomology calculations in this subsection can be done with the
aid of Theorem~\ref{thm:line_bundle_cohomology} below; in essence, one
has to check that certain subcomplexes of the sphere~$S^{n-1}$ are
contractible. This is what is behind the calculations in the paper
\cite{H-nonlin_toric} which, however, uses a dual point of view, using
the fact that the fans considered above are normal fans of
polytopes. We omit the details.

\subsubsection*{Ishida's formula}

If the support function~$h$ is upper convex (equivalently, if the
associated line bundle is generated by global sections), the negatives
of the linear functions $h_\sigma$ for $n$-dimensional cones $\sigma
\in \fan$ span a polytope~$Q$ in~$M_\bR$ with vertices in~$M$. Since
\[H^k (X_{\fan_P}; L_h) = \cases{0 & if \( k \not= 0\)\ , \cr
\bigoplus_{M \cap Q} \bC & if \(k=0\)\ ,}\]
Theorem~\ref{thm:generalised_Brion} specialises to
\cite[Theorem~2.3]{Ishida-polyhedral} for complete fans.

\subsubsection*{An explicit example}

\begin{example}
\label{ex:explicit}
We consider the case $n=2$, $N = \bZ^2$ and $N_\bR = \bR^2$. Let
$\fan$ be the unique complete fan in $N_\bR$ 
whose $1$-cones are generated by the following four vectors:
\[v_1 = \pmatrix{1 \cr 1} \qquad v_2 = \pmatrix{0 \cr 1} \qquad v_3 =
\pmatrix{-1 \cr 1} \qquad v_4 = \pmatrix{0 \cr -1}\]
Let $X, Y \in M = \hom_\bZ (N, \bZ)$ denote the dual of the standard
basis of $N = \bZ^2$.
Let $h \colon N_\bR \rTo \bR$ be the support function specified by the values
\[h(v_1)= 0 \qquad h(v_2)= -2 \qquad h(v_3)= 0 \qquad h(v_4)= -2\]
given by extending linearly over cones; for example, on $\sigma = \mathrm{cone} (v_1,
  v_2)$ it agrees with the linear function $h_\sigma = 2X-2Y \in M$ which corresponds to a
\textsc{Laurent\/} monomial written $x^2 y^{-2}$ . Using
Theorem~\ref{thm:line_bundle_cohomology} we can
explicitly compute the cohomology of $L_h$ (see end of
  \S\ref{sec:cohom-line-bundle} below). It turns out that $H^0 (X_\fan;
  L_h) = 0$, and that $\dim H^2 (X_\fan; L_h) = 1$ concentrated in homogeneous
degree~$-Y \in M$. The space $H^1 (X_\fan; L_h)$ is $4$-dimensional, with a
  $1$-dimensional contribution coming from degrees $0$, $-X+Y$,
  $Y$ and~$X+Y$. The right-hand side of Equation~(\ref{eq:main_res})
  thus is (denoting the indeterminates again by~$x$ and~$y$)
  \[-1-x^{-1}y - y -xy + y^{-1} \ .\]
The left-hand side is worked out easily as well. For example, the summand
corresponding to $\sigma = \mathrm{cone} (v_1, v_2)$ is the rational function
represented by the lattice point enumerator of the shifted cone
\[-h_\sigma + \sigma^\vee = (-2X+2Y) + \mathrm{cone} \left(\, X, -X+Y
  \,\right) \ \subset M_\bR\]
or, in coordinates of $M_\bR \iso \bR^2$,
\[-h_\sigma + \sigma^\vee = \pmatrix{-2 \cr 2} + \mathrm{cone} \left(\, \pmatrix{1 \cr 0}, \pmatrix{-1
    \cr 1} \, \right) \ .\]
This lattice point enumerator is given by
\[{{x^{-2} y^2} \over {(1-x^{-1}y)(1-x)}} \ .\]
In total, the left-hand side of Equation~(\ref{eq:main_res}) equals
\[{{x^{-2} y^2} \over {(1-x^{-1}y)(1-x)}} +
{{x^2 y^2} \over {(1-x^{-1})(1-xy)}} \qquad\qquad\qquad \relax\]
\[\qquad\qquad\qquad + {{x^{-2} y^{-2}} \over
  {(1-x^{-1})(1-x^{-1}y^{-1})}} + {{x^2 y^{-2}} \over {(1-x)(1-xy^{-1})}}\]
\smallbreak
\noindent which coincides with the \textsc{Laurent\/} polynomial $-1-x^{-1}y -
y -xy + y^{-1}$, as an explicit calculation shows.
\end{example}

\section{Cohomology of line bundles}
\label{sec:cohom-line-bundle}

\subsubsection*{\v Cech cohomology of quasi-coherent sheaves}

\begin{numpar}
  \label{par:cell-structure}
  Let $\fan$ be a finite complete fan in~$N_\bR$. By taking
  intersection of positive-dimensional cones with the unit sphere
  $S^{n-1}$ (defined with respect to any inner product) the fan
  induces the structure of a regular $CW$ complex on~$S^{n-1}$. Given
  a cone $\sigma \in \fan$ we write $\bar \sigma = \sigma \cap
  S^{n-1}$ for the corresponding cell of~$S^{n-1}$; this includes the
  case of the empty cell $\bar 0$. We fix once and for all
  orientations of the cells and write $[\bar \sigma : \bar \tau]$ for
  the incidence number of~$\bar\sigma$ and~$\bar\tau$. By convention,
  we have $[\bar \tau : \bar 0] = 1$ for all $1$-dimensional cones
  $\tau \in \fan$. Regularity of the $CW$ decomposition implies that
  $[\bar \sigma : \bar \tau] \in \{-1,0,1\}$ for all cones $\sigma,
  \tau \in \fan$. Note that in the (augmented) cellular chain complex
  of~$S^{n-1}$ the empty cell corresponds to the augmentation
  (concentrated in degree $-1$).
\end{numpar}

The computation of \textsc{\v Cech\/} cohomology does not depend on using the
complex numbers as a ground field. So let $A$ denote a commutative ring, and
let $X_\fan$ denote the toric $A$-scheme associated to~$\fan$; it is obtained
by gluing the affine $A$-schemes $\mathrm{Spec}\, A[M \cap \sigma^\vee]$ where
$\sigma$ varies over the elements of~$\fan$.

\begin{proposition}
  \label{prop:CechCohomology}
  Let $\mathcal{F}$ be a quasi-coherent sheaf on~$X_\fan$. Then we can
  compute the cohomology modules $H^k (X_\fan; \mathcal{F})$ as the
  cohomology of the \textsc{\v Cech} cochain complex
  $\mathcal{C^\bullet} =(\mathcal{C}^\bullet, d)$ which is defined by
  \[\mathcal{C}^d = \bigoplus_{{\sigma \in \fan} \atop {\mathrm{codim}\, \sigma = d}}
  \mathcal{F}^\sigma\ ,\]
  with differential defined on direct summands by
  \[\mathcal{F}^\sigma \rTo[l>=3em]^{[\bar\sigma:\bar\tau]}
  \mathcal{F}^\tau \ .\]
  Here $\mathcal{F}^\sigma$ denotes the module of sections
  of~$\mathcal{F}$ over the affine open subset of~$X_\fan$ determined
  by the cone~$\sigma$. In particular, $\mathcal{F}^\sigma$ is an
  $A[M \cap \sigma^\vee]$-module. The family
  $(\mathcal{F}^\sigma)_{\sigma \in \fan}$ of modules
  determines $\mathcal{F}$ completely.
\end{proposition}

This variant of \textsc{\v Cech} cohomology should be well-known, but
unfortunately there seems to be no published proof available. For the
reader's convenience we give a proof in \S\ref{sec:proof-cech} below.

\subsubsection*{Torus-equivariant line bundles}

\begin{numpar}
  \label{par:support_function}
  We will apply Proposition~\ref{prop:CechCohomology} in the case
  where $\mathcal{F}$ is a torus equivariant algebraic line bundle
  on~$X_\fan$.  Recall \cite[\S2]{Oda-Toric} that such a sheaf is
  specified by a support function $h \colon N_\bR \rTo \bR$ which is
  linear on each cone, and takes integral values on~$N$; in other
  words, for each $\sigma \in \fan$ there exists $h_\sigma \in
  \hom_\bZ (N, \bZ)=M$ such that $h|_\sigma = \left(h_\sigma \tensor
    \mathrm{id}_\bR\right) |_\sigma$.  The linear function $h_\sigma$
  is well-defined up to addition of a linear function which vanishes
  on $N \cap \sigma$.
\end{numpar}

\begin{numpar}
  \label{par:L_h}
  The line bundle $L_h$ corresponding to a support function~$h$ has a
  very explicit description: On the affine open set corresponding to
  $\sigma \in \fan$ the space of sections is the free $A[M \cap
  \sigma^\vee]$-module of rank~$1$ with basis $-h_\sigma$. Note that
  all these modules are contained in the free $A$-module $A[M]$ with
  basis $M$, hence we may consider them as $M$-graded $A$-modules.
\end{numpar}

\begin{numpar}
  \label{par:hom_piece_cohomology}
  We can apply Proposition~\ref{prop:CechCohomology} to the line bundle~$L_h$.
  The resulting cochain complex of $A$-modules has a natural
  $M$-grading, and the differentials are homogeneous of degree~$0$ with
  respect to this grading (in the language of free modules, all the terms in
  $\mathcal{C}^\bullet$ have a basis consisting of a subset of~$M$, and all
  structure maps are induced by inclusion of subsets.) Hence the cohomology
  modules $H^k (X_\fan; L_h)$ have a direct sum decomposition
  \[H^k (X_\fan; L_h) = \bigoplus_{\mathbf{b} \in M} H^k (X_\fan;
  L_h)_\mathbf{b} \ ,\]
  with $H^k (X_\fan; L_h)_\mathbf{b}$ being isomorphic to the
  cohomology of the \hbox{degree-$\mathbf{b}$} sub-cochain complex
  $\mathcal{C}_\mathbf{b}^\bullet = (\mathcal{C}_\mathbf{b}^\bullet, d)$
  of~$\mathcal{C}^\bullet$.
\end{numpar}

\begin{numpar}
  \label{par:desc_hom_chain_cx}
  The cochain complex $\mathcal{C}_\mathbf{b}^\bullet$ itself admits a simple
  description: It is given by
  \[\mathcal{C}_\mathbf{b}^d = \bigoplus_{{\sigma \in \fan,\ \mathrm{codim}\, \sigma = d}
    \atop {\mathbf{b} + h_\sigma \in \sigma^\vee}} A \]
  with differential induced by incidence numbers as before. Now if
  $\tau$ is a face of $\sigma \in \fan$ then $\sigma^\vee \subseteq
  \tau^\vee$, so $\mathbf{b}+h_\sigma \in \sigma^\vee$ implies $\mathbf{b}+h_\sigma \in
  \tau^\vee$, hence $\mathbf{b}+h_\tau \in \tau^\vee$ (for $h_\tau - h_\sigma
  \in \tau^\vee$, and $\tau^\vee$ is closed under addition since it is
  a convex cone). Thus the set
  \begin{equation}
    \label{eq:def_S_h_u}
    S(h,\mathbf{b}) := \bigcup_{{\sigma \in \fan} \atop
      {\mathbf{b}+h_\sigma \in \sigma^\vee}} \bar \sigma
  \end{equation}
  is a sub-complex of~$S^{n-1}$, and $\mathcal{C}_\mathbf{b}^\bullet$
  is nothing but the augmented cellular chain complex
  of~$S(h,\mathbf{b})$, re-indexed suitably as a cochain complex. In
  other words, we have shown:
\end{numpar}

\begin{theorem}
  \label{thm:line_bundle_cohomology}
  Suppose $\fan$ is a complete fan in~$N_\bR$, and $h \colon N_\bR
  \rTo \bR$ is a support function on~$\fan$. Let $L_h$ denote the
  algebraic line bundle on~$X_\fan$ associated to~$h$ (\ref{par:L_h}),
  and define the space $S(h,\mathbf{b})$ as in~{\rm
    (\ref{eq:def_S_h_u})}. For all $\mathbf{b} \in M$ there is an
  isomorphism of $A$-modules
  \[H^k(X_\fan; L_h)_\mathbf{b} \iso \tilde H_{n-1-k}
  (S(h,\mathbf{b}); A)\]
  where $\tilde H_d (\,\cdot\,; A)$ denotes reduced cellular (or singular)
  homology with coefficients in~$A$.
  \qed
\end{theorem}

For this to make sense, it is imperative to consider the augmented
cellular chain complex to compute $\tilde H_d$ with augmentation
concentrated in degree~$-1$. In other words, $\tilde H_{-1}
(\emptyset) = A$ by convention, while $\tilde H_{-1} (X) = 0$
whenever $X \not= \emptyset$.

\medbreak

The advantage of Theorem~\ref{thm:line_bundle_cohomology} over the standard
result as given in \cite[Theorem~2.6]{Oda-Toric} is that the former deals with
the cell complex $S(h,\mathbf{b})$ arising as the intersection of a
sub-fan of~$\fan$ with $S^{n-1}$, whereas the latter relies on
computing certain subsets of~$N_\bR$ with a rather more delicate
combinatorial structure.

\medbreak

The theorem leads immediately to some general observations. For example, the
remark following Theorem~\ref{thm:line_bundle_cohomology} implies:

\begin{corollary}
  The top-dimensional cohomology is given by
  \[H^n (X_\fan; L_h)_\mathbf{b} = \cases{0 & if there exists \(\sigma \in
    \fan,\ \sigma \not= \{0\} \) with \(\mathbf{b} + h_\sigma \in \sigma^\vee\) \cr A
    & otherwise.} \]
  Moreover, if $H^n (X_\fan; L_h)_\mathbf{b} = A$, then $H^k
  (X_\fan; L_h)_\mathbf{b} = 0$ for all $k \not= n$. \qed
\end{corollary}

Suppose now that $K$ is a subcomplex of $S^{n-1}$. Then $\tilde
H_{n-1} (K; A) \not= 0$ if and only if $K = S^{n-1}$. Indeed, if $K
\not= S^{n-1}$ then $K$ misses an $(n-1)$-dimensional cell
of~$S^{n-1}$, \ie, there exists an $n$-dimensional cone $\sigma \in
\fan$ such that $K$ is contained in $S^{n-1} \setminus \mathrm{int}\,
\bar \sigma$. Now $S^{n-1} \setminus \mathrm{int}\, \bar \sigma$ is
contractible, hence has trivial reduced homology. The homology long
exact sequence of the pair $(K, S^{n-1} \setminus \mathrm{int}\, \bar
\sigma)$ proves the assertion.

If there exists $\mathbf{b} \in M$ such that
$\mathcal{S}(h,\mathbf{b}) = S^{n-1}$, then $\mathbf{b}$ is contained
in the intersection of the closed half-spaces $-h_\rho + \rho^\vee$
where $\rho$ varies over the $1$-dimensional cones in~$\fan$. Since
$\fan$ is complete, this implies that for all $\mathbf{a} \in \bZ^n$
there exists $\rho \in \fan$ with $\mathbf{a} \in -h_\rho +
\rho^\vee$, thus $\mathcal{S}(h,\mathbf{a}) \not= \emptyset$.
Conversely, if $\mathcal{S}(h,\mathbf{b}) = \emptyset$ for some
$\mathbf{b} \in M$ then there is no $\mathbf{a} \in M$ with
$\mathcal{S}(h,\mathbf{a}) = S^{n-1}$ (in fact, there is a
$1$-dimensional cone $\rho \in \fan$ with $\mathbf{a} \notin -h_\rho +
\rho^\vee$). Together with Theorem~\ref{thm:line_bundle_cohomology}
this shows that the line bundle $L_h$ cannot have global sections and
$n$th cohomology and the same time:

\begin{corollary}
  At least one of the $A$-modules $H^0 (X_\fan; L_h)$ and $H^n (X_\fan;
  L_h)$ is trivial. \qed
\end{corollary}

\subsubsection*{On Example~\ref{ex:explicit}}

Recall the notation from Example~\ref{ex:explicit}; we will use the
field $A = \bC$ of complex numbers.  To work out the complex
$\mathcal{S}(h,\mathbf{b}) \subseteq S^1$ for given $\mathbf{b} \in M$
one can start from a sketch of the halfspace arrangement $-h_{\rho_j}
+ \rho_j^\vee$, $j=1,2,3,4$ given by the shifted duals of the
$1$-dimensional cones in~$\fan$.  Furthermore, it is enough to
consider those~$\mathbf{b}$ which are contained in some bounded region
of the resulting decomposition of~$M_\bR$ since $H^* (X_\fan; L_h)$ is
finite-dimensional.

In our example, this leaves us to check contributions from five
elements of~$M$ only. We will use coordinate notation for this
paragraph. It is easily verified that $\mathcal{S} (h, (0,-1)^t) =
\emptyset$, so $H^2 (X_\fan; L_h)_{(0,-1)^t} = \tilde H_{-1}
(\emptyset) = \bC$. If $\mathbf{b}$ is one of the vectors $\{(0,0)^t,\
(-1,1)^t,\ (0,1)^t,\ (1,1)^t\}$, then $\mathcal{S}(h,\mathbf{b})$ is a
$0$-sphere corresponding to the intersection of the cones spanned by
$v_1$ and~$v_3$ with the unit sphere in~$N_\bR = \bR^2$. Thus $H^1
(X_\fan; L_h)_\mathbf{b} = \tilde H_0 (S^0) = \bC$ in these cases.

\section{Proof of Theorem~\ref{thm:generalised_Brion}}

The proof of Theorem~\ref{thm:generalised_Brion} proceeds by verifying
a \textsc{Laurent} series identity first. Let as before $h \colon
N_\bR \rTo \bR$ be a support function, and choose corresponding linear
functions $h_\sigma \in M$ for $\sigma \in \fan$ (\ref{par:support_function}).
Define a formal \textsc{Laurent} power series
\begin{equation}
  \label{eq:def_R_fan_series}
  R[\fan, h] = \sum_{\sigma \in \fan} (-1)^{\mathrm{codim}\, \sigma}
  R[-h_\sigma+\sigma^\vee] \ .
\end{equation}
Fix $\mathbf{a} \in \bZ$; we want to consider the coefficient
of~$x^\mathbf{a}$ in~$R[\fan, h]$. The summand corresponding to
$\sigma \in \fan$ contributes $0$ if $\mathbf{a} + h_\sigma \notin
\sigma^\vee$, and it contributes $(-1)^{\mathrm{codim}\, \sigma}$
otherwise. Since $0^\vee = \bR^n$ we get a contribution of~$(-1)^n$
for $\sigma = 0$ always. In other words, the coefficient
of~$x^\mathbf{a}$ is the \textsc{Euler} characteristic of the chain
complex~$\mathcal{C}^\bullet_\mathbf{a}$
(\ref{par:desc_hom_chain_cx}), using $A = \bC$ again:
\[\sum_{k=0}^n (-1)^k \dim_\bC\, \mathcal{C}^k_\mathbf{a} = \chi
(\mathcal{C}^\bullet_\mathbf{a}) \]
The \textsc{Euler} characteristic can be computed using the cohomology
groups of the cochain complex as well, so the coefficient of
$x^\mathbf{a}$ is given by
\begin{equation}
  \label{eq:Euler_char_hom}
  \chi (\mathcal{C}^\bullet_\mathbf{a}) = \sum_{k=0}^n (-1)^k
  \dim_\bC\, \tilde H_{n-1-k} S(h,\mathbf{a}) \ .
\end{equation}
Using Theorem~\ref{thm:line_bundle_cohomology} we see that this is
equal to $\sum_{k=0}^n (-1)^k \dim_\bC\, H^k (X_\fan;
L_h)_\mathbf{a}$. Since the cohomology of $L_h$ is finitely generated
(the variety $X_\fan$ is complete by hypothesis), we see that this
coefficient is zero for almost all $\mathbf{a} \in \bZ^n$. In
particular, $R[\fan, h]$ is a \textsc{Laurent} polynomial.

\medskip

Let $\Pi$ denote the $P$-submodule of $S$ generated by the rational functions
corresponding to rational polyhedral cones. According to
\cite[Theorem~1.2]{Ishida-polyhedral} there is a unique $P$-linear
homomorphism (here $Q(P)$ denotes the quotient field of~$P$ as before)
\[\rho \colon \Pi \rTo Q(P)\]
with $\rho (R[b+\sigma]) = R(b+\sigma)$ for all $b \in M_\bR$ and all
pointed rational polyhedral cones $\sigma \subset M_\bR$ (see also
\cite[Theorem~2.4]{Beck:Brion_formulae}). Note that $\rho$ preserves
\textsc{Laurent} polynomials as they are finite sums of
\textsc{Laurent} power series associated to sets of the form
\hbox{$\mathbf{a} + \{0\}$}. In particular, $\rho (x^\mathbf{b}) =
x^\mathbf{b} \in P \subset Q(P)$ for all $\mathbf{b} \in M$.---If the
rational polyhedral cone $\sigma$ contains a line, 
then it can be shown that $\rho (K[\sigma]) = 0$,
cf.~\cite[Lemma~2.1]{Ishida-polyhedral} or
\cite[Lemma~2.5]{Beck:Brion_formulae}.

\medskip

We now apply the homomorphism~$\rho$ to the \textsc{Laurent} power
series $R[\fan, h]$. On the one hand, we have
\begin{eqnarray*}
  \rho (R[\fan, h]) &=& \sum_{\sigma \in \fan} (-1)^{\mathrm{codim}\, \sigma}
  \rho (R[-h_\sigma+\sigma^\vee]) \\
  &=& \sum_{{\sigma \in \fan} \atop {\dim \sigma =
      n}} (-1)^{\mathrm{codim}\, \sigma} \rho (R[-h_\sigma+\sigma^\vee]) \\
  &=& \sum_{{\sigma \in \fan} \atop {\dim \sigma = n}}
  R(-h_\sigma+\sigma^\vee) \\
  &=& R(\fan, h) \ .
\end{eqnarray*}
(The second equality comes from the fact that if $\mathrm{codim}\,
\sigma >0$, then the dual cone $\sigma^\vee$ contains a line.)
On the other hand, we have seen that $R[\fan, h]$ is a
\textsc{Laurent} polynomial. Hence $R(\fan, h) = \rho(R[\fan, h]) = R[\fan, h]$ is
a \textsc{Laurent} polynomial as well, and as seen before the coefficient of
$x^\mathbf{a}$ is given by $\chi (H^*(X_\fan; L_h)_\mathbf{a})$. This
finishes the proof.

\bigbreak

As a final remark, we can also use Equation~(\ref{eq:Euler_char_hom})
to identify the coefficients of the monomials in~$R(\fan, h)$ as this
is an intermediate step in the above proof. The result then reads:

\begin{corollary}
  The coefficient of $x^\mathbf{a}$ in $R(\fan, h)$ is equal to
  $(-1)^{n-1} \tilde \chi (S(h,\mathbf{a}))$, the reduced
  \textsc{Euler} characteristic of the cell complex $S(h,
  \mathbf{a})$. In other words,
  \[R(\fan, h) = (-1)^{n-1} \sum_{\mathbf{a} \in \bZ^n} \tilde \chi
  (S(h,\mathbf{a})) \cdot x^\mathbf{a} \ .\qed\]
\end{corollary}

\section{Proof of Proposition~\ref{prop:CechCohomology}}
\label{sec:proof-cech}

Let $A$ denote a commutative ring with unit. For a
complete fan $\fan$ in~$N_\bR$ we let $X_\fan$ denote the associated
toric scheme defined over~$A$. A quasi-coherent sheaf of modules
$\mathcal{F}$ on~$X_\fan$ determines, by evaluation on affine pieces,
a diagram of $A$-modules
\[D(\mathcal{F}) = D \colon \fan^\mathrm{op} \rTo A\mathrm{-Mod}, \quad \sigma
\mapsto D^\sigma = \mathcal{F}^\sigma\]
(where as before $\mathcal{F}^\sigma$ denotes the $A$-module of
sections of~$\mathcal{F}$ over the open affine subset of~$X_\fan$
determined by~$\sigma$, cf.~Proposition~\ref{prop:CechCohomology}).
The functor $\mathcal{F} \mapsto D(\mathcal{F})$ is exact: A short
exact sequence of quasi-coherent sheaves yields a short exact sequence
of diagrams. We need the fact that we can compute sheaf cohomolgy by
higher derived limits of the associated diagram:

\begin{lemma}
  There are canonical isomorphisms
  \[H^j (X_\fan; \mathcal{F}) \iso \lim_\leftarrow \!^j D(\mathcal{F})
  \ .\]
\end{lemma}

\begin{proof}
  Given a cone $\sigma \in \fan$ write $U_\sigma$ for the open affine
  subset of~$X_\fan$ determined by~$\sigma$. Then by construction
  $D(\mathcal{F})^\sigma = \Gamma (U_\sigma; \mathcal{F})$, and the
  case $j=0$ of the Lemma is just the sheaf axiom: A global section is
  uniquely determined by a collection of compatible local section.

  Recall now that sheaf cohomology
  can be computed with flasque resolutions. That is, considering
  $\mathcal{F}$ as a sheaf of \textsc{abel}ian groups, choose a
  resolution
  \begin{equation}
    \label{eq:resolution}
    \mathcal{F} \rTo \mathcal{G}_0 \rTo \mathcal{G}_1 \rTo \cdots
  \end{equation}
  with all the $\mathcal{G}_i$ being flasque. Let $U \subseteq X_\fan$
  be an open subset; then $\mathcal{G}_i |_{U}$ is flasque, and $H^j
  (U; \mathcal{F})$ is isomorphic to the cohomology groups of the
  cochain complex
  \begin{equation}
    \label{eq:compute_sheaf_cohomology}
    \Gamma (U; \mathcal{G}_0|_U) \rTo \Gamma (U; \mathcal{G}_1|_U)
  \rTo \cdots \ .
  \end{equation}
  Passing to associated diagrams of \textsc{abel}ian groups the
  resolution~(\ref{eq:resolution}) gives rise to a cochain complex
  \begin{equation}
    \label{eq:resolution_diagram}
    D(\mathcal{F}) \rTo D(\mathcal{G}_0) \rTo D(\mathcal{G}_1) \rTo
  \cdots \ .
  \end{equation}
  We claim that this is in fact a resolution of~$D(\mathcal{F})$,
  considered as a diagram of \textsc{abel}ian groups. Indeed, given a
  cone $\sigma \in \fan$ the cochain complex
  \[D(\mathcal{G}_0)^\sigma \rTo D(\mathcal{G}_1)^\sigma \rTo \cdots\]
  is nothing but the cochain
  complex~(\ref{eq:compute_sheaf_cohomology}) for $U = U_\sigma$,
  hence its $j$th cohomology group is isomorphic to $H^j (U_\sigma:
  \mathcal{F})$. But $U_\sigma$ is affine and $\mathcal{F}$
  quasi-coherent, so these groups vanish for $j \geq 1$, proving the
  claim.

  We observe that the resolution~(\ref{eq:resolution_diagram}) is
  flasque in the sense that the canonical restriction maps
  \begin{equation}
    \label{eq:can_restriction}
    D(\mathcal{G}_i)^\sigma \rTo \lim_{\tau \subset \sigma}
    D(\mathcal{G}_i)^\tau
  \end{equation}
  are surjective. Indeed, using the definition of associated diagrams,
  the map~(\ref{eq:can_restriction}) corresponds to the restriction
  map
  \[\Gamma (U_\sigma; \mathcal{G}_0) \rTo \Gamma (\bigcup_{\tau
    \subset \sigma} U_\tau; \mathcal{G}_i)\]
  which is surjective since $\mathcal{G}_i$ is flasque. Hence we can
  use the resolution~(\ref{eq:resolution_diagram}) to compute higher
  derived inverse limuits of~$D(\mathcal{F})$ by applying the functor
  $\lim$ to~(\ref{eq:resolution_diagram}), then taking cohomology
  groups. However, applying $\lim$ to~(\ref{eq:resolution_diagram})
  yields precisely the cochain
  complex~(\ref{eq:compute_sheaf_cohomology}) for $U = X_\fan$, which
  computes $H^j (X_\fan; \mathcal{F})$. Taking into account the
  well-known fact that the higher derived inverse limits of a diagram
  of $A$-modules can be computed in the category of diagrams of
  \textsc{abel}ian groups, we have thus proved the Lemma.
\end{proof}

\bigbreak

The proof of Proposition~\ref{prop:CechCohomology} thus reduces to
proving the following claim:

\begin{proposition}
  \label{prop:alternate_cech}
  Let $D \colon \fan^\mathrm{op}\rTo A\mathrm{-Mod}$, $\sigma \mapsto
  D^\sigma$ be a diagram of $A$-modules, where $A$ is an arbitrary ring with
  unit. Form the cochain complex $\mathcal{C}^\bullet =
  \mathcal{C}(D)^\bullet$ by setting
  \[\mathcal{C}^k = \mathcal{C}(D)^k = \bigoplus_{{\sigma \in \fan} \atop {\mathrm{codim}\, \sigma
  = k}} D^\sigma \ ,\]
  with differential defined on direct summands by
  \[D^\sigma \rTo[l>=4em]^{[\bar \sigma : \bar \tau]} D^\tau \ .\]
  Then the \textsc{\v Cech} cohomology modules $\check H^k (D) = h^k
  (\mathcal{C}(D)^\bullet)$ are naturally isomorphic to the higher derived
  inverse limits $\lim\!^k (D)$.
\end{proposition}

\medskip

The proof of the proposition will occupy the rest of this section.
First, for $n=1$ we know that $\fan$ consists
of the zero-cone, $\bR_{\geq 0}$ and $\bR_{\leq 0}$. The cochain
complex $\mathcal{C}^\bullet$ has the form
\[D^{\bR_{\geq 0}} \oplus D^{\bR_{\leq 0}} \rTo^+ D^{\{0\}} \ ,\]
and the result is well-known in this case.

We can thus restrict to the case $n \geq 2$. We extend the cell
structure on $S^{n-1}$ introduced in~\ref{par:cell-structure} to a regular
cell structure on $B^n$ with a single $n$-cell denoted~$\bar B$.
The canonical maps $\lim (D) \rTo D^\sigma$, modified by the incidence numbers
$[\bar B : \bar\sigma]$, assemble to a map
\[\iota \colon \lim (D) \rTo \bigoplus_{{\sigma \in \fan}
  \atop {\dim \sigma = n}} D^\sigma = \mathcal{C}^0\]
which is, by the properties of incidence numbers, a co-augmentation of
the cochain complex~$\mathcal{C}^\bullet$.

\begin{lemma}
  \label{lem:iota_iso}
  The map $\iota$ is injective and induces an isomorphism $\lim (D)
  \iso \check H^0 (D)$.
\end{lemma}

\begin{proof}
  An element of $\lim (D)$ is determined by its images in
  the~$D^\sigma$ where $\sigma$ ranges over all $n$-dimensional cones
  of~$\fan$. Conversely, an element of~$\mathcal{C}^0$ lies in the kernel of
  $\mathcal{C}^0 \rTo \mathcal{C}^1$ if and only if its components in $D^\sigma$ and
  $D^\tau$ agree in $D^{\sigma \cap \tau}$ where $\sigma, \tau \in
  \fan$ are $n$-dimensional cones with $(n-1)$-dimensional
  intersection. Such an element thus determines a unique element
  of~$\lim D$ mapping to the given element of~$\mathcal{C}^0$.
\end{proof}

\medbreak

Observe now that the functor $D \mapsto \check H^* (D)$ is a
$\delta$-functor \cite[\S2.1]{Weibel-Intro}. Indeed, a short exact
sequence of diagrams
\begin{equation}
  \label{eq:ses_sheaves}
  0 \rTo D \rTo E \rTo F \rTo 0
\end{equation}
gives rise to a short exact sequence of cochain complexes, hence by
the snake lemma to an associated natural long exact sequence in
cohomology. Since $D \mapsto \lim\!^* (D)$ is a universal
$\delta$-functor \cite[\S2.1 and \S2.5]{Weibel-Intro}, it follows that
we have uniquely determined natural maps
\[\nu_k \colon \lim\!^k (D) \rTo \check H^k (D)\]
such that $\nu_0$ is the isomorphism of Lemma~\ref{lem:iota_iso}, and such
that the $\nu_k$ give rise to a commutative ladder diagram in cohomology for
every short exact sequence of the form~(\ref{eq:ses_sheaves}).

\medbreak

To prove that the maps $\nu_k$ are isomorphisms, we consider a
decreasing filtration of the diagram~$D$. For $0 \leq j \leq n$ we
write
\[\kappa_j D \colon \fan^\mathrm{op} \rTo  A\mathrm{-Mod}, \quad
\sigma \mapsto \cases{D^\sigma & if \(\mathrm{codim}\, \sigma \leq j\),
  \cr 0 & else.}\]

\begin{lemma}
  \label{lem:iso_tau_0}
  The maps $\nu_k$ are isomorphisms for all diagrams of the form $\kappa_0
  D$.
\end{lemma}

\begin{proof}
  The diagram $\kappa_0 D$ has non-zero values only on $n$-dimensional
  cones, hence $\check H^k (\kappa_0 D) = 0$ for $k > 0$.  It is easy to
  check that $\lim\!^k (\kappa_0 D) = 0$ for $k>0$ (for example, examine
  the cochain complex of \cite[Vista~3.5.12]{Weibel-Intro} which
  computes higher derived inverse limits). Thus $\nu_k \colon
  \lim\!^k (\kappa_0 D) \rTo \check H^k (\kappa_0 D)$ is an isomorphism
  for all~$k$.
\end{proof}

\medbreak

We proceed by induction on~$j$ and state the {\bf induction
  hypothesis:\/} {\it The maps $\nu_k \colon \lim\!^k (\kappa_{j-1} D) \rTo
\check H^k (\kappa_{j-1} D)$ are isomorphisms for all $k \geq 0$ and
all diagrams~$D$.} The case $j=1$ is covered by the previous Lemma.

\medbreak

We have a sequence of epimorphisms of diagrams
\[D = \kappa_n D \rTo^{e_n} \kappa_{n-1} D \rTo^{e_{n-1}} \ldots
\rTo^{e_1} \kappa_0 D\]
and consequently a collection of short exact sequences (for $1 \leq j
\leq n$)
\[0 \rTo \ker (e_j) \rTo \kappa_j D \rTo \kappa_{j-1} D \rTo 0 \ .\]
Consider the associated ladder diagram for some fixed $k \geq 1$:
{\small
\begin{diagram}
  \lim\!^{k-1} (\kappa_{j-1} D) & \rTo & \lim\!^k (\ker e_j) & \rTo &
  \lim\!^k (\kappa_j D) & \rTo & \lim\!^k (\kappa_{j-1} D) & \rTo &
  \lim\!^{k+1} (\ker e_j) \\
  \dTo>{\nu_{k-1}} && \dTo>{\nu_k} && \dTo>{\nu_k} && \dTo>{\nu_k} &&
  \dTo>{\nu_{k+1}} \\
  \check H^{k-1} (\kappa_{j-1} D) & \rTo & \check H^k (\ker e_j) & \rTo &
  \check H^k (\kappa_j D) & \rTo & \check H^k (\kappa_{j-1} D) & \rTo &
  \check H^{k+1} (\ker e_j)
\end{diagram}
}

By our induction hypothesis we know that the first and fourth
vertical arrow are isomorphisms.  In view of the $5$-lemma it is
enough to show that the second and fifth vertical arrow are
isomorphisms as well. Since $\tau_n D = D$ this proves the assertion
of Proposition~\ref{prop:alternate_cech}.

\medbreak

We are left to show that the maps $\nu_k \colon \lim\!^k (\ker e_j)
\rTo \check H^k (\ker e_j)$ are isomorphisms for all~$k$. Now the
diagram $\ker e_j$ has non-trivial entries only on cones of
codimension~$j$, and can thus be written as a direct sum of ``atomic''
diagrams with a single non-trivial entry. Since both $\lim\!^k$ and
$\check H^k$ commute with direct sums of atomic diagrams (the former
for abstract reasons, the latter by direct inspection), the induction
step is completed if we can verify the following assertion:

\begin{lemma}
 \label{lem:atomic_C_iso}
  Let $C$ be an $A$-module.
  Let $\tau \in \fan$ be a cone of codimension $j > 0$, and let
  $C_\tau$ denote the atomic $\fan^\mathrm{op}$-diagram with
  non-trivial value~$C$ attained at~$\tau$. Then the maps $\nu_k
  \colon \lim\!^k (C_\tau) \rTo \check H^k (C_\tau)$ are
  isomorphisms for all $k \geq 0$.
\end{lemma}

(Note that by uniqueness of the natural maps~$\nu_k$ the direct sum
decomposition of the diagram $\ker e_j$ carries over to a direct sum
decomposition of the corresponding~$\nu_k$.)

\medbreak

The Lemma will follow from a brute-force calculation.  By
construction, $\check H^k (C_\tau) = 0$ for $k \not=j$, and $\check
H^j (C_\tau) = C$. We have $\lim (C_\tau) = 0$ since $\tau$ has positive
codimension. To compute the higher derived inverse limits, we embed
the diagram $C_\tau$ into a short exact sequence
\begin{equation}
  \label{eq:ses_c_tau}
  0 \rTo C_\tau \rTo C_{\geq \tau} \rTo C_{>\tau} \rTo 0
\end{equation}
where we write
\[C_{\geq \tau} \colon \fan^\mathrm{op} \rTo A\mathrm{-Mod},
\quad \sigma \mapsto \cases{C & if \(\sigma \supseteq \tau\)\ , \cr
  0 & else}\]
(with non-trivial structure maps identities), and the diagram
$C_{>\tau}$ is defined similarly.

\begin{lemma}
  \label{lem:restrict_same}
  The higher derived inverse limits remain unchanged when the diagram
  $C_{\geq\tau}$ is restricted to the subcategory $\tau \downarrow
  \fan = \{ \sigma \in \fan \,|\, \sigma \supseteq \tau\}$. More
  precisely, the canonical restriction maps
  $\lim_{\fan^\mathrm{op}}\!^k\, C_{\geq\tau} \rTo \lim_{(\tau
    \downarrow \fan)^\mathrm{op}}\!\!\!\!^k\, (C_{\geq\tau}|_{(\tau
    \downarrow \fan)})$
  are isomorphisms.
\end{lemma}

\begin{proof}
  This can be read off from the usual cochain complex computing higher
  derived inverse limits as given in
  \cite[Vista~3.5.12]{Weibel-Intro}.
\end{proof}

\begin{lemma}
  \label{lem:const_ball}
  We have $\lim (C_{\geq \tau}) = C$ and $\lim\!^k
  (C_{\geq \tau}) = 0$ for $k \geq 1$.
\end{lemma}

\begin{proof}
  By Lemma~\ref{lem:restrict_same} we get an isomorphism $\lim\!^k
  (C_{\geq \tau}) \iso H^k (N (\tau\downarrow\fan)^\mathrm{op}; C)$ for
  all $k \geq 0$, where $N$ denotes the nerve of the category. But
  $\tau \downarrow \fan$ has an initial object, hence is contractible.
\end{proof}

\begin{lemma}
  \label{lem_h_C_tau}
  We have $\lim\!^k (C_\tau) = 0$ for $k \not=j$, and
  $\lim\!^j (C_\tau) = C$.
\end{lemma}

\begin{proof}
  This follows from the long exact sequence associated to the short
  exact sequence~(\ref{eq:ses_c_tau}) and the calculations in the
  previous two lemmas.
\end{proof}

As a consequence it is enough to consider the case $k=j$ in
Lemma~\ref{lem:atomic_C_iso}.

\begin{lemma}
  \label{lem_Hcheck_C_tau}
  We have $\check H^k (C_{\geq \tau}) = 0$ for $k>0$.
\end{lemma}

\begin{proof}
  In short, this follows from the fact that the cochain complex
  $\mathcal{C}(C_{\geq \tau})^\bullet$ is a re-indexed variant of the
  cellular chain complex computing the reduced homology of a
  $(j-1)$-sphere with coefficients in~$C$, so $\check H^k (C_{\geq
    \tau}) \iso \tilde H_{j-1-k} (S^{j-1}; C)$.

  In more detail, the poset $\tau \downarrow \fan$ is known to be
  isomorphic to the $j$-dimensional quotient fan $\fan/\tau$ as
  defined in \cite[Corollary~1.7]{Oda-Toric}. It is a complete fan in
  the vector space $N_\bR / \mathrm{span}(\tau)$, its cones are given
  by the images of $\sigma \in \tau \downarrow \fan$ under the
  quotient map $N_\bR \rTo N_\bR / \mathrm{span}(\tau)$. The fan
  $\fan/\tau$ induces a cell structure on some unit sphere $S^{j-1}$
  in $N_\bR / \mathrm{span}(\tau)$, and taking the incidence numbers
  coming from   the fan~$\fan$ as defined before, we see that
  $\mathcal{C}(C_{\geq\tau})^\bullet$ is, up to re-indexing, an
  augmented cellular chain complex of~$S^{j-1}$. This chain complex is
  slightly non-standard: The augmentation maps are given by
  $\id_C$ or $-\id_C$, depending on the incidence numbers $[\bar\sigma
  : \bar\tau]$. However, it is not difficult to show that this chain
  complex is isomorphic to a standard chain complex for any choice of
  orientations of the cones in~$\fan/\tau$, the required isomorphism
  being constructed by induction of the dimension of the cones,
  starting with $\tau$. We omit the details.
\end{proof}

\medbreak

We are now ready to prove Lemma~\ref{lem:atomic_C_iso}. Consider the
following piece of the ladder diagram relating $\lim\!^k$ and $\check
H^*$:
\begin{equation}
  \label{eq:final_diagram}
  \begin{diagram}
    \lim\!^{j-1} (C_{\geq \tau}) & \rTo &
    \lim\!^{j-1} (C_{>\tau}) & \rTo[l>=3em]^f & \lim\!^j (C_\tau) & \rTo &
    \lim\!^j (C_{\geq \tau}) = 0 \\
    \dTo<\iso>{\nu_{j-1}} &&
    \dTo<\iso>{\nu_{j-1}} && \dTo>{\nu_j} \\
    \check H^{j-1} (C_{\geq \tau}) & \rTo &
    \check H^{j-1} (C_{>\tau}) & \rTo^g & \check H^j (C_\tau) &
    \rTo & \check H^j (C_{\geq \tau}) = 0
  \end{diagram}
\end{equation}
Both rows are exact. The entries on the right are trivial by
Lemma~\ref{lem:const_ball} and direct inspection of the cochain
complex $\mathcal{C} (C_{\geq \tau})^\bullet$, respectively. The first
vertical map is an isomorphism; for $j=1$ it is the map $\nu_0$, and
for $j>1$ it follows from Lemmas~\ref{lem:const_ball}
and~\ref{lem_Hcheck_C_tau}. The second vertical map is an isomorphism
in view of our induction hypothesis (note that $C_{>\tau} =
\kappa_{j-1} C_{\geq \tau}$). From the Five Lemma we conclude that the
third vertical map is an isomorphism as desired.

{\small
\raggedright
\bibliographystyle{alpha}

\begin{thebibliography}{Oda88}

\bibitem[BHS]{Beck:Brion_formulae}
Matthias Beck, Christian Haase, and Frank Sottile.
\newblock {Theorems of Brion, Lawrence, and Varchenko on rational generating
  functions for cones}.
\newblock arXiv:math.CO/0506466.

\bibitem[Bri88]{Brion-latticepoints}
Michel Brion.
\newblock Points entiers dans les poly\`edres convexes.
\newblock {\em Ann. Sci. \'Ecole Norm. Sup. (4)}, 21(4):653--663, 1988.

\bibitem[BV97]{Brion-Vergne}
Michel Brion and Mich{\`e}le Vergne.
\newblock Lattice points in simple polytopes.
\newblock {\em J. Amer. Math. Soc.}, 10(2):371--392, 1997.

\bibitem[H{\"u}ta]{H-nonlin_toric}
Thomas H{\"u}ttemann.
\newblock {{$K$}-Theory of non-linear projective toric varieties}.
\newblock arXiv:math.KT/0508431.

\bibitem[H{\"u}tb]{H-Brion}
Thomas H{\"u}ttemann.
\newblock On a theorem of {B}rion.
\newblock arXiv:math.CO/0607297.

\bibitem[Ish90]{Ishida-polyhedral}
Masa-Nori Ishida.
\newblock Polyhedral {L}aurent series and {B}rion's equalities.
\newblock {\em Internat. J. Math.}, 1(3):251--265, 1990.

\bibitem[Oda88]{Oda-Toric}
Tadao Oda.
\newblock {\em Convex bodies and algebraic geometry}, volume~15 of {\em
  Ergebnisse der Mathematik und ihrer Grenzgebiete (3) [Results in Mathematics
  and Related Areas (3)]}.
\newblock Springer-Verlag, Berlin, 1988.
\newblock An introduction to the theory of toric varieties, Translated from the
  Japanese.

\bibitem[Wei94]{Weibel-Intro}
Charles~A. Weibel.
\newblock {\em An introduction to homological algebra}, volume~38 of {\em
  Cambridge Studies in Advanced Mathematics}.
\newblock Cambridge University Press, Cambridge, 1994.

\end{thebibliography}

}

\end{document}